\documentclass[10pt]{amsart}

\usepackage{epsfig}
\usepackage{graphicx}
\usepackage{graphics}
\usepackage{psfrag}
\usepackage{amsmath}
\usepackage{amsthm}
\usepackage{amsfonts}
\usepackage{amssymb}
\usepackage{amscd}
\usepackage[all]{xy}

\usepackage{color}
\usepackage{hyperref}

\newcommand{\comment}[1]{}

\newcommand {\hide}[1]{}
\theoremstyle{theorem}
\newtheorem{theorem}{Theorem}[section]
\newtheorem{corollary}[theorem]{Corollary}
\newtheorem{lemma}[theorem]{Lemma}
\newtheorem{proposition}[theorem]{Proposition}

\theoremstyle{definition}

\newtheorem{definition}[theorem]{Definition}
\newtheorem{example}[theorem]{Example}

\newtheorem{remark}[theorem]{Remark}

\newcommand {\C}     {\mathbb{C}}
\newcommand {\cH}     {\mathcal{H}}
\newcommand {\cP}     {\mathcal{P}}
\newcommand {\cQ}     {\mathcal{Q}}

\newcommand {\eps}     {\varepsilon}

\newcommand {\PP}     {\mathbb{P}}

\newcommand {\R}     {\mathbb{R}}
\newcommand {\Z}     {\mathbb{Z}}
\newcommand {\Zer}     {\mbox{\rm Zer}}
\newcommand {\Sphere}{\mbox{${\bf S}$}}     
\newcommand {\Ball}{\mbox{${\bf B}$}}     

\begin{document}

\title[A sharper estimate on the Betti numbers]
{A sharper estimate on the Betti numbers of sets defined by 
quadratic inequalities}

\author{Saugata Basu}
\author{Michael Kettner}
\address{School of Mathematics,
Georgia Institute of Technology, Atlanta, GA 30332, U.S.A.}
\thanks{The first author was supported in part by an NSF Career Award 0133597 and
a Alfred P. Sloan Foundation Fellowship. The second author was partially 
supported by the European RTNetwork 
Real Algebraic and Analytic Geometry, Contract No.~HPRN-CT-2001-00271.}
\email{\{saugata,mkettner\}@math.gatech.edu}


\date{\today}



\begin{abstract}
In this paper we consider the problem of bounding the Betti numbers,
$b_i(S)$, of
a semi-algebraic set $S \subset \R^k$ defined by polynomial inequalities
$P_1 \geq 0,\ldots,P_s \geq 0$, where $P_i \in \R[X_1,\ldots,X_k]$ and
$\deg(P_i) \leq 2$, for $1 \leq i \leq s$. We prove that for $0\le i\le k-1$, 
\[
b_i(S) \le\frac{1}{2}\left(
\sum_{j=0}^{min\{s,k-i\}}{{s}\choose j}{{k+1}\choose {j}}2^{j}\right).
\]
In particular, for $2\le s\le \frac{k}{2}$, we have
\[
b_i(S)\le \frac{1}{2} 3^{s}{{k+1}\choose {s}} 
\leq \frac{1}{2} \left(\frac{3e(k+1)}{s}\right)^s.
\]
This improves the bound of $k^{O(s)}$ proved by Barvinok in \cite{Barvinok}.
This improvement is made possible by a new approach, whereby
we first bound the Betti numbers of 
non-singular complete intersections of complex projective varieties
defined by generic quadratic forms,
and use this bound to obtain bounds in the real semi-algebraic case.
\end{abstract}

\maketitle
\section{Introduction}
%
The topological complexity of semi-algebraic sets, measured by 
their Betti numbers (ranks of their  singular homology groups), 
has been the subject of many investigations. 
For any topological space $X$, we will denote by $b_i(X) = b_i(X,\Z_2)$
the $i$-th Betti number of $X$ with $\Z_2$-coefficients, and we will denote
by $b(X)$ the sum $\sum_{i \geq 0} b_i(X)$. 
Note that, since the homology groups of a semi-algebraic
set $S \subset \R^k$ are finitely generated, it follows from the 
Universal Coefficients Theorem,
that
$b_i(S,\Z_2) \geq b_i(S,\Z)$, where $b_i(S,\Z)$ are the ordinary Betti numbers
of $S$ with integer coefficients (see \cite{Hatcher}, Corollary 3.A6 (b)).
Hence, the bounds proved in this
paper also apply  to the ordinary Betti numbers. 
(The use of $\Z_2$ 
coefficients is necessitated by our use of Smith inequalities 
in the proof of the main theorem.)

The initial result on bounding the Betti numbers of semi-algebraic sets 
defined by polynomial inequalities was proved independently by 
Oleinik and Petrovskii~\cite{OP}, Thom~\cite{Thom} and Milnor~\cite{Milnor}.
They proved:
\begin{theorem}\cite{OP,Thom,Milnor}
\label{the:OP} 
Let 
\[
\cP=\{P_1,\ldots,P_s\}\subset\R[X_1,\ldots,X_k]
\] 
with 
$\deg(P_i)\le d$, $1\le i\le s$ and let $S\subset\R^k$ be the set 
defined by 
\[
P_1\geq 0,\ldots,P_s\geq 0.
\]
Then, 
\[
b(S)=O(sd)^k.
\]
\end{theorem}
Notice that the above bound is exponential in $k$ and this exponential 
dependence is unavoidable (see Example \ref{eg:example} below). See also 
\cite{B03,GV05,BPR05} for more recent work extending the above bound
to more general classes of semi-algebraic sets.
%
\subsection{Semi-algebraic Sets Defined by Quadratic Inequalities}
In this paper we consider a restricted class of semi-algebraic sets - namely, 
semi-algebraic sets defined by quadratic inequalities. 
Since sets defined by linear 
inequalities have no interesting topology, sets defined by quadratic 
inequalities can be considered to be the simplest class of semi-algebraic 
sets which can have  non-trivial topology. Such sets are in fact quite 
general, as every semi-algebraic set can be defined by (quantified) 
formulas involving  only quadratic polynomials (at 
the cost of increasing the number of variables and the size of the formula). 
Moreover, as in the case of general semi-algebraic sets, the Betti numbers 
of such sets can be exponentially large as can be seen in the following 
example.

\begin{example}
\label{eg:example}
The set~$S\subset\R^k$ defined by 
\[
X_1(X_1-1)\ge0,\ldots, X_k(X_k-1)\ge0,
\]
has $b_0(S)=2^k$.
\end{example}

However, it turns out that for a semi-algebraic set $S \subset \R^k$
defined by $s$ quadratic inequalities,
it is possible to obtain upper bounds on the Betti numbers of $S$ 
which are polynomial in $k$ and exponential only in $s$.
The first such result was proved by Barvinok who proved the following 
theorem.

\begin{theorem}\cite{Barvinok}
\label{the:barvinok}
Let $S \subset \R^k$ be defined by $P_1 \geq 0, \ldots,  P_s \geq 0$, 
$\deg(P_i) \leq 2, 1 \leq i \leq s$. Then,
$b(S) \leq k^{O(s)}$.
\end{theorem}

Theorem~\ref{the:barvinok} is proved using a duality argument 
that interchanges the roles of $k$ and $s$, 
and reduces the original problem to that of bounding the Betti numbers of
a semi-algebraic set in $\R^s$ defined by $k^{O(1)}$ polynomials of degree 
at most $k$. One can then use Theorem~\ref{the:OP} to obtain a bound
of $k^{O(s)}$. The constant hidden in the exponent of the above bound
is at least two. Also, the bound in Theorem~\ref{the:barvinok} is polynomial
in $k$ but exponential in $s$. 
The exponential dependence on $s$
is unavoidable as remarked in \cite{Barvinok}, but the implied constant
(which is at least two) in the exponent of Barvinok's bound is not optimal.

Using Barvinok's result, as well as inequalities derived from the 
Mayer-Vietoris sequence, a polynomial bound (polynomial both in
$k$ and $s$) was proved in \cite{B03} on the top few Betti numbers of
a set defined by quadratic inequalities. More precisely the
following theorem is proved there.

\begin{theorem}\cite{B03}
\label{the:quadratic}
Let $\ell > 0$ and ${\rm R}$ a real closed field. 
Let  $S \subset {{\rm R}}^k$ be defined by 
\[
P_1 \ge 0,\ldots, P_s \ge 0,
\] with
$\deg(P_i) \leq 2$.
Then, \[
b_{k-\ell}(S) \leq {s \choose {\ell}} k^{O(\ell)}.
\]
\end{theorem}

Notice that for fixed $\ell$, the bound in Theorem \ref{the:quadratic} is
polynomial in both $s$ as well as $k$.

Apart from their intrinsic mathematical interest (in distinguishing the
semi-algebraic sets defined by quadratic inequalities from general
semi-algebraic sets), 
the bounds in Theorems \ref{the:barvinok} and \ref{the:quadratic} have
motivated recent work on designing polynomial time algorithms for computing
topological invariants of semi-algebraic sets defined by quadratic
inequalities (see \cite{Barvinok2,GP,B06a,B06b,BZ}). Traditionally
an important goal in algorithmic semi-algebraic geometry has been to design
algorithms for computing topological invariants
of semi-algebraic sets, 
whose worst-case complexity matches the best upper bounds known for
the quantity being computed.  
It is thus of interest to tighten the bounds on the Betti numbers
of semi-algebraic sets defined by quadratic inequalities,
as has been done recently 
in the case of general semi-algebraic sets (see for example
\cite{GV05,B03,BPR05}).
%
\subsection{Brief Outline of Our Method}
In this paper we use a new method to bound the Betti numbers of 
semi-algebraic sets defined by quadratic inequalities. 
Our method is to first bound the Betti numbers of 
complex projective varieties which are 
non-singular complete intersections 
defined by quadratic forms in general position. 
It is a well known fact from complex geometry, 
(see for instance, \cite{Lewis}, pp. 122)
that the Betti numbers of a complex projective variety which is a 
non-singular complete intersection
depend only on the sequence of degrees of the polynomials defining the
variety.
Moreover, there exist precise formulas for the Betti numbers of 
such varieties, using well-known techniques
from algebraic geometry  (see Theorem \ref{the:betti} below).

Our strategy for bounding the Betti numbers of semi-algebraic sets
in $\R^k$ defined by $s$ quadratic inequalities is as follows. Using certain
{\em infinitesimal deformations} we first reduce the problem to bounding
the Betti numbers of another closed and {\em bounded}  
semi-algebraic set defined by a new family quadratic polynomials.  
We then  use inequalities obtained from the
Mayer-Vietoris exact sequence to further reduce the problem of bounding
the Betti numbers of this new  semi-algebraic set, 
to the problem of bounding the
Betti numbers of the {\em real projective varieties} defined by
each
$\ell$-tuple, $\ell \leq s$, of the new polynomials.  
The new family of polynomials also has the
property that the {\em complex projective variety}  defined by each
$\ell$-tuple, $\ell \leq k$, of these polynomials is a {\em non-singular
complete intersection}.
As mentioned above we have precise information about the Betti numbers
of these complex complete intersections. An application of {\em 
Smith inequalities} then allows us to obtain bounds on the Betti 
numbers of the real parts
of these varieties and as a result on the Betti numbers of the
original semi-algebraic set.
Because of the direct nature of our proof we are able to remove the constant
in the exponent in the bounds proved in \cite{Barvinok,B03} and this 
constitutes the main contribution of this paper.

\begin{remark}
We remark here that the technique used in this paper  
was proposed as a possible alternative method 
by Barvinok in \cite{Barvinok}, who did not pursue this 
further in that paper. 
Also,
Benedetti, Loeser, and Risler \cite{BLR91}
used a similar technique for  proving 
upper bounds on the number of connected components 
of real algebraic sets in $\R^k$
defined by polynomials of degrees bounded by $d$. However, these
bounds (unlike the situation considered in this paper) 
are exponential in $k$.
Finally, there exists another possible  method for 
bounding the Betti numbers of semi-algebraic
sets defined by quadratic inequalities, using a spectral sequence argument
due to Agrachev \cite{Agrachev}. However, this  
method also produces a non-optimal bound
of the form $k^{O(s)}$ (similar to Barvinok's bound)
where the constant in the exponent is at least two. We omit the
details of this argument referring the reader to 
\cite{B06a} for an indication of the proof (where the case of
computing, and as a result, bounding the Euler-Poincar\'e characteristics
of such sets is worked out in full details).
\end{remark}

We prove the following theorem.
\begin{theorem}\label{the:Pgre0}
Let $\cP=\{P_1,\ldots,P_s\}\subset\R[X_1,\ldots,X_k]$, 
$s\le k$. 
Let $S\subset \R^k$ be defined by
\[
P_1\ge0,\ldots,P_s\ge0
\]
with $\deg(P_i)\le 2$. Then, for $0\le i\le k-1$, 
\[
b_i(S) \le\frac{1}{2}\left(
\sum_{j=0}^{min\{s,k-i\}}{{s}\choose j}{{k+1}\choose {j}}2^{j}\right).
\]
In particular, for $2\le s\le \frac{k}{2}$, we have
\[
b_i(S)\le \frac{1}{2} 3^{s}{{k+1}\choose {s}} 
\leq \frac{1}{2} \left(\frac{3e(k+1)}{s}\right)^s.
\]
\end{theorem}
As a consequence of the proof of Theorem~\ref{the:Pgre0} we get a new bound 
on the sum of the Betti numbers, which we state for the sake of completeness.
\begin{corollary}
\label{cor:bsum}
Let $\cP=\{P_1,\ldots,P_s\}\subset\R[X_1,\ldots,X_k]$, 
$s\le k$. 
Let $S\subset \R^k$ be defined by
\[
P_1\ge0,\ldots,P_s\ge0
\]
with $\deg(P_i)\le 2$. Then, 
\[
b(S) \le \frac{1}{2}k\left(
\sum_{j=0}^{s}{{s}\choose j}{{k+1}\choose {j}}2^{j}\right).
\]
\end{corollary}

The rest of the paper is organized as follows. In Section~\ref{sec:prelim},
we recall some well known results from complex algebraic geometry on the
Betti numbers of non-singular complex projective varieties which
are complete intersections, and also some classical results from algebraic
topology which we need for the proof of our main result. In
Section~\ref{sec:proof}, we prove Theorem~\ref{the:Pgre0}. 
Finally, in Section \ref{sec:open} we state some open problems.
\section{Mathematical Preliminaries} \label{sec:prelim}
%
In this section we recall a few basic facts about Betti numbers and 
complex projective varieties which are non-singular 
complete intersections, as well as fix some notations.

Throughout the paper, $\PP_{\R}^k$ (respectively, $\PP_{\C}^k$) 
denotes the real (respectively, complex) projective space of dimension~$k$, 
$\Sphere^k_r$ (resp., $\Ball^{k+1}_r$) denotes the sphere (resp., 
closed ball)  centered at the origin and of radius $r$ in $\R^{k+1}$.
For any polynomial $P \in \R[X_1,\ldots,X_k]$,
we denote by $P^h \in \R[X_1,\ldots,X_{k+1}]$ the homogenization of $P$ 
with respect to $X_{k+1}$.

For any family of polynomials
$\cP=\{P_1,\ldots,P_{\ell}\}\subset\R[X_1,\ldots,X_k]$, and 
$Z \subset \R^k$,
we denote by $\Zer(\cP,Z)$ 
the set of common zeros of $\cP$ in $Z$. Moreover, 
for any family of homogeneous polynomials 
$\cQ=\{Q_1,\ldots,Q_{\ell}\}\subset\R[X_1,\ldots,X_{k+1}]$, 
we denote by $\Zer(\cQ,\PP_{\R}^k)$ (resp., $\Zer(\cQ,\PP_{\C}^k)$)
the set of common zeros of $\cQ$ in $\PP_{\R}^k$ 
(resp., $\PP_{\C}^k$).
%
%
\subsection{Some Results from Algebraic Topology}
%
Let 
\[
\mathcal{Q}=\{Q_1,\ldots,Q_{\ell}\}\subset\R[X_1,\ldots,X_{k+1}]
\] 
be a set of  
homogeneous polynomials. Denote by $\mathcal{Q}_J$ the 
set $\{Q_j| j\in J\}$ for $J\subset\{1,\ldots,\ell\}$.
We have the following inequality which is a consequence of the
Mayer-Vietoris exact sequence.

\begin{proposition}\label{prop:spec}
Let $Z\subset\R^{k+1}$. For $0\le i\le k-1$,
\[
b_i(\bigcup_{j=1}^{\ell}\Zer(Q_j,Z))\le
\sum_{j=1}^{i+1}\sum_{J\subset\{1,\ldots,\ell\},|J|=j} 
b_{i-j+1}(\Zer(\mathcal{Q}_J,Z)).
\]
\end{proposition}
\begin{proof}
See \cite{B03}, Lemma~2. 
\end{proof}

We also use the well-known Alexander duality theorem which relates
the Betti numbers of a compact subset of a sphere to those of its
complement.

\begin{theorem} [Alexander Duality]\label{the:alex}
Let $r>0$. For any closed subset~$A\subset S_r^k$,
\[
H_i(S_r^k\setminus A)\approx \tilde{H}^{k-i-1}(A),
\]
where $\tilde{H}^i(A)$, $0\le i\le k-1$, denotes the 
reduced cohomology group of $A$.
\end{theorem}
\begin{proof}
See \cite{Massey}, Theorem 6.6.
\end{proof}
Finally, we state a version of the 
Smith inequality which plays a crucial role in the proof of the main
theorem. 
Recall that for any compact topological 
space equipped with an involution, inequalities derived from the {\em Smith
exact sequences} allows one to bound the {\em sum} of the Betti numbers 
(with $\Z_2$ coefficients) 
of the fixed point set of the involution by the sum of the Betti numbers 
(again with $\Z_2$ coefficients) 
of the space itself (see for instance, \cite{Viro}, pp. 131). 
In particular, we have for a complex projective
variety defined by real forms, with the involution taken to be 
complex conjugation, the following theorem.
\begin{theorem}[Smith inequality]\label{the:smith}
Let ${\mathcal Q} \subset \R[X_1,\ldots,X_{k+1}]$ be a family of
homogeneous polynomials.
Then,
\[
b(\Zer({\mathcal Q},\PP^k_{\R}))\le b(\Zer({\mathcal Q},\PP^k_{\C})).
\]
\end{theorem}

%
\subsection{Complete Intersection Varieties}\label{ssec:complint}
%
\begin{definition}
A projective variety~$X$ of codimension~$n$ is a 
\textit{non-singular complete intersection} 
if it is the intersection of $n$ non-singular 
hypersurfaces that meet transversally at each point of the 
intersection.
\end{definition}
Fix  an $j$ tuple of natural numbers $\bar{d} = (d_1,\ldots,d_j)$. Let
$X_{\C} = \Zer(\{Q_1,\ldots,Q_j\},\mathbb{P}_{\C}^{k})$, 
such that the degree of $Q_i$ is $d_i$, 
denote a complex projective variety of  
codimension~$j$ which is a non-singular complete intersection.

Let $b(j,k,\bar{d})$ denote the sum of the Betti numbers with 
$\mathbb{Z}_2$ coefficients of $X_{\C}$. This is well defined since the 
Betti numbers only depend only on the degree sequence 
and not on the specific $X_{\C}$.

The function $b(j,k,\bar{d})$ satisfies the following (see \cite{BLR91}):

\[
b(j,k,\bar{d}) = 
\begin{cases}
c(j,k,\bar{d}) & \mbox{ if } k-j \mbox{ is even,} \\
2(k -j +1) - c(j,k,\bar{d}) & \mbox{ if } k-j \mbox{ is odd},
\end{cases}
\]
where
\[
c(j,k,\bar d) = 
\begin{cases} 
k + 1 & \mbox{ if } j = 0, \\
d_1\ldots d_j & \mbox{ if } j = k, \\
d_k c(j-1,k-1,(d_1,\ldots,d_{k-1})) - (d_k-1)c(j,k-1,\bar{d}) & \mbox{ if } 
0 < j < k. 
\end{cases}
\]

In the special case when each $d_i = 2$, we denote by
$b(j,k)=b(j,k,(2,\ldots,2))$. 
We then have the following recurrence
for $b(j,k)$.

\[
b(j,k) = 
\begin{cases} 
q(j,k) & \mbox{ if } k-j \mbox{ is even}, \\
2(k-j+1)-q(j,k) &\mbox{ if } k-j \mbox{ is odd}, 
\end{cases}
\]
where 
\[
q(j,k) = 
\begin{cases} 
k + 1 & \mbox{ if } j = 0, \\
2^j & \mbox{ if } j = k, \\
2 q(j-1,k-1) - q(j,k-1) & \mbox{ if } 0 < j < k. 
\end{cases}
\]
Next, we show some properties of $q(j,k)$.
\begin{lemma}
\label{lemma:prop_q}
\item
\begin{enumerate}
\item $q(1,k)= k + 1/2(1-(-1)^k)$ and $q(2,k)=(-1)^k k+ k$.
\item For $2\le j \le k$, $|q(j,k)|\le 2^{j-1}{k\choose{j-1}}$.
\item For $2\le j \le k$ and $k-j$ odd, $2(k-j+1)-q(j,k)\le 2^{j-1}{k\choose{j-1}}$.
\end{enumerate}
\end{lemma}
\begin{proof}
The first part is shown by two easy computations and note that 
$2(k-2+1)-q(2,k)=2k-2$ if $k-2$ is odd. 
Hence, we can assume that the statements are true for $k-1$ and that $3\le j < k$. 
Note that for the special case $j=k-1$, we have that 
$2^{k-1}\le 2^{k-2} {{k-1}\choose {k-2}}$ since 
$k > 2$.Then, 
\begin{eqnarray*}
|q(j,k)| & = & |2q(j-1,k-1) - q(j,k-1)|\\
 & \le & 2|q(j-1,k-1)| + |q(j,k-1)|\\
 & \le & 2\cdot 2^{j-2}{{k-1}\choose{j-2}} + 2^{j-1}{{k-1}\choose j-1}\\
 & = & 2^{j-1}{k\choose{j-1}}.
\end{eqnarray*}
and, for $k-j$ odd, 
\begin{eqnarray*}
2(k-j+1)-q(j,k) & = & 2(k-j+1) - 2q(j-1,k-1) + q(j,k-1)\\
 & \le & |2((k-1)-(j-1)+1) - q(j-1,k-1)|\\
 &      & + |q(j-1,k-1)| + |q(j,k-1)|\\
 & \le & 2^{j-2}{{k-1}\choose {j-2}} + 2^{j-2}{{k-1}\choose {j-2}} + 2^{j-1}{{k-1}\choose {j-1}}\\
 & \le & 2^{j-1}\left({{k-1}\choose {j-2}}+{{k-1}\choose {j-1}}\right) = 2^{j-1}{k\choose {j-1}}.
\end{eqnarray*}
\end{proof}
\pagebreak[2]
Hence, we get the following bound for $b(j,k)$.
\begin{theorem}
\label{the:betti}
\begin{enumerate}
\item[]
\item 
$
b(1,k) = 
\begin{cases} 
q(0,k-1) & \mbox{ if } k \mbox{ is even}, \\
q(0,k) &\mbox{ if } k \mbox{ is odd}, 
\end{cases}
$
\item $b(j,k)\le 2^{j-1}{k\choose{j-1}}$, for $2\le j \le k$.
\end{enumerate}
\end{theorem}
\begin{proof}
Follows from Lemma~\ref{lemma:prop_q}.
\end{proof}
\begin{proposition}
\label{prop:general}
Let $\eps>0$ and let $P_{\eps}=(\frac{2}{\eps})^2 - \sum_{i=1}^{k+1} X_i^2$. Then, 
there exist a 
family~${\cH_{\eps}=\{H_{\eps,1},\ldots,H_{\eps,s}\}}\subset\R[X_1,\ldots,X_{k+2}]$, 
$s\le k$, of positive definite quadratic forms such that   
$\Zer(\cH_{\eps,J}\cup\{P^h_{\eps}\},\PP^{k+1}_{\C})$ 
is a non-singular complete intersection 
for every $J\subset\{1,\ldots,s\}$.
\end{proposition}

\begin{proof}
Note that for any family~$\cH=\{H_1,\ldots,H_s\}\subset\R[X_1,\ldots,X_{k+2}]$ 
of quadratic forms such that their coefficients 
are algebraically independent over 
$\mathbb{Q}$, 
$\Zer(\cH_J,\PP^{k+1}_{\C})$, 
$J\subset\{1,\ldots,s\}$, 
is a non-singular complete intersection by Bertini's Theorem 
(see \cite{Harris}, Theorem~17.16). 
Moreover, recall that the set of positive definite quadratic forms 
is open in the set of quadratic forms (over $\R$). 
Thus, we can choose for every~$\eps>0$ a 
family~$\cH_{\eps}=\{H_{\eps,1},\ldots,H_{\eps,s}\}\subset\R[X_1,\ldots,X_{k+2}]$, 
$s\le k$, of positive definite quadratic forms such that their coefficients 
are algebraically independent over 
$\mathbb{Q}(\eps)$, 
and   
$\Zer(\cH_{\eps,J}\cup\{P^h_{\eps}\},\PP^{k+1}_{\C})$ 
will be a non-singular complete intersection 
for every ${J\subset\{1,\ldots,s\}}$.
\end{proof}
The following proposition allows to replace a set of real quadratic
forms by another family obtained by infinitesimal perturbations
of the original family and whose zero sets are non-singular complete
intersections in complex projective space.
\begin{proposition}\label{prop:compl}
Let $\eps>0$ and let 
\[
P_{\eps}=(\frac{2}{\eps})^2 - \sum_{i=1}^{k+1} X_i^2.
\] 
Let 
\[
\cQ=\{Q_1,\ldots,Q_{s}\}\subset\R[X_1,\ldots,X_{k+2}],
\] 
$s\le k$, be a set of quadratic forms and 
let 
\[
{\cH_{\eps}=\{H_{\eps,1},\ldots,H_{\eps,s}\}}\subset\R[X_1,\ldots,X_{k+2}]
\] 
be 
a family of positive definite quadratic forms such that   
$\Zer(\cH_{\eps,J}\cup\{P^h_{\eps}\},\PP^{k+1}_{\C})$ 
is a non-singular complete intersection  for every $J\subset\{1,\ldots,s\}$.
 
For $t \in \C$,
let 
$$
\displaylines{
\tilde{\cQ}_{\eps,t}= \{\tilde{Q}_{\eps,t,1},\ldots,\tilde{Q}_{\eps,t,s}\}
\;\mbox{with} \cr
\tilde{Q}_{\eps,t,i} = (1-t)Q_i + t H_{\eps,i}.
}
$$
Then, for all sufficiently small 
$\delta > 0$, and any $J \subset \{1,\ldots,s\}$,
\[
\Zer(\tilde{\cQ}_{\eps,\delta,J}\cup\{P^h_{\eps}\},\PP_{\C}^{k+1})
\]
is a non-singular complete intersection.
\end{proposition}
\begin{proof}
Let $J \subset \{1,\ldots,s\}$, and
let $T_J \subset \C$ be defined by,
$$
\displaylines{
T_J = \{ t \in \C \;\mid \; \Zer(\tilde{\cQ}_{\eps,t,J},\PP_{\C}^{k+1})\;\mbox{is a
non-singular complete intersection} \; \}. 
}
$$
Clearly, $T_J$ contains $1$. Moreover, since being a non-singular complete 
intersection is stable condition, 
$T_J$ must contain an open neighborhood of $1$ in $\C$ and so must
$T = \cap_{J\subset \{1,\ldots,s\}} T_J$.
Finally, the set $T$ is constructible,
since it can be defined by a first order formula. 
Since a constructible subset of $\C$
is either finite  or the complement of a finite set
(see for instance, \cite{BPR05}, Corollary 1.25),
$T$ must contain an interval $(0,t_0), t_0 >0$. 
\end{proof}
%
\section{Proof of Theorem~\ref{the:Pgre0}}
\label{sec:proof}
%
In order to prove Theorem~\ref{the:Pgre0} we need what follows next:

Let $\mathcal{P}=\{P_1,\ldots,P_s\}\subset\R[X_1,\ldots,X_k]$, $s\le k$, with 
$\deg(P_i)\le 2$, $1\le i\le s$. Let 
$S\subset\R^k$ be the basic semi-algebraic set 
defined by ${P_1\ge0,\ldots,P_s\ge 0}$, 
and let 
\[
S_{\eps}=S \cap \Ball^k_{1/\eps}.
\]
\begin{proposition}\label{prop:SvstildeS}
For all sufficiently small $\eps >0$,
the homology groups of $S$ and $S_{\eps}$ are isomorphic. 
Moreover, $S_{\eps}$ is bounded.
\end{proposition}
\begin{proof}
The proof follows from Hardt's triviality theorem (see \cite{BPR06}, 
Theorem~5.45.) 
and is similar to the proof of Lemma~1 in \cite{B99}.
\end{proof}
\begin{figure}
\begin{center}
\includegraphics[scale=0.4]{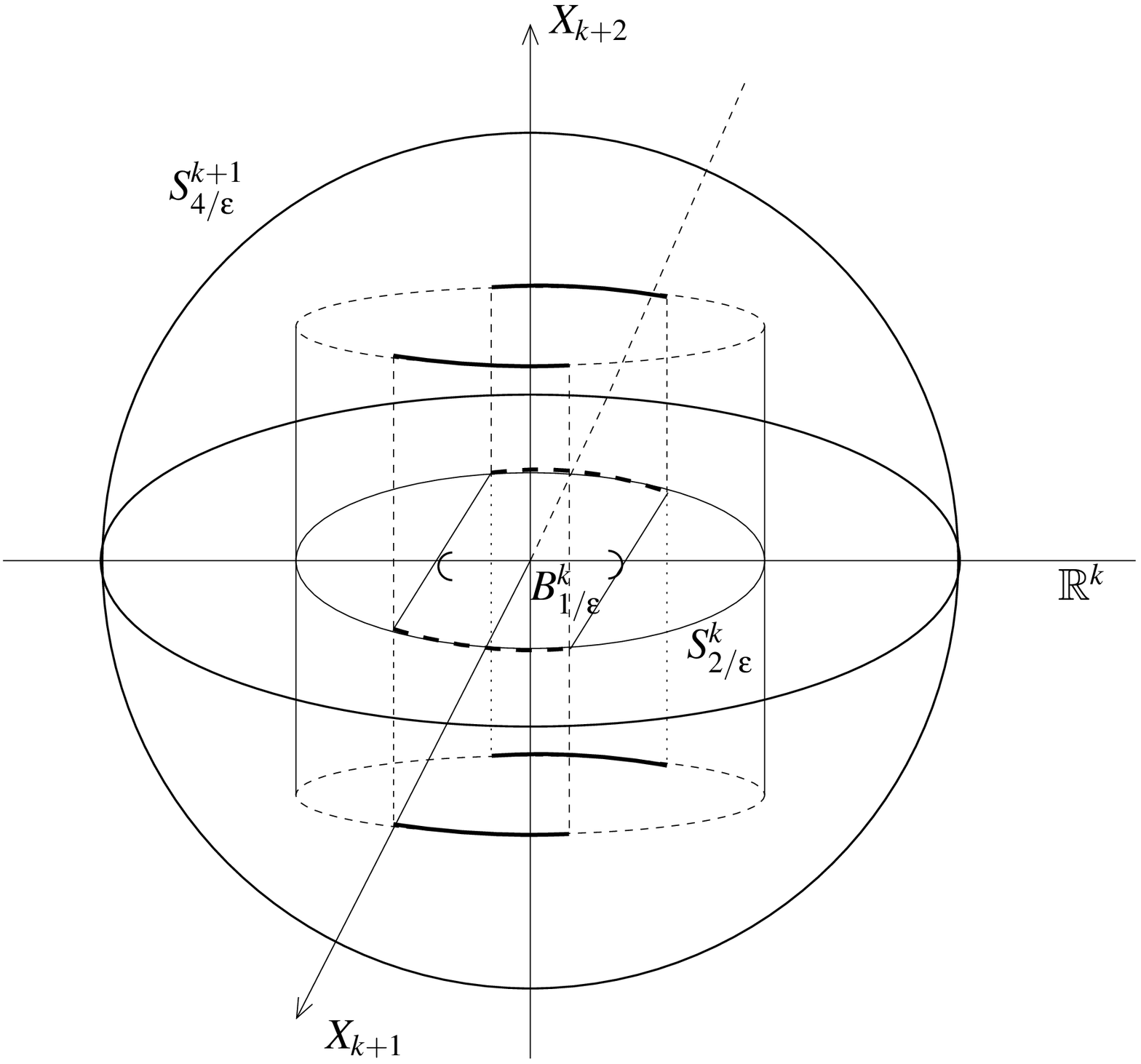} 
\end{center}
\caption{lifting the ball~$\Ball^k_{1/\eps}$ onto the sphere~$\Sphere^{k+1}_{4/\eps}$}
\label{fig:lemma}
\end{figure}
Before we continue, consider Figure~\ref{fig:lemma} which will be 
helpful for the 
following. The cylinder $\Ball^k_{1/\eps}\times\R$ above the 
ball $\Ball^k_{1/\eps}\subset\R^k$ intersects the sphere~$\Sphere^k_{2/\eps}$ in two disjoint 
copies (each homeomorphic to $\Ball^k_{1/\eps}$). 
Each cylinder 
above those copies 
intersects the sphere~$\Sphere^{k+1}_{4/\eps}$ in two disjoint copies. 
Thus, there are four disjoint copies of $\Ball^k_{1/\eps}$ on the 
sphere~$\Sphere^{k+1}_{4/\eps}$ 
(each homeomorphic to $\Ball^k_{1/\eps}$). Notice
that each such copy does not intersect the 
equator of the sphere $\Sphere^{k+1}_{4/\eps}$ in $\R^{k+2}$ (i.e.
the set $\Sphere^{k+1}_{4/\eps}\cap\Zer(X_{k+2},\R^{k+2})$).

Let $S^h_{\eps}$ be the basic semi-algebraic set defined by 
$P^h_1\ge 0,\ldots, P^h_s\ge 0$ contained in 
$C_{\eps}$, where 
\[
C_{\eps}={(\Ball^k_{1/\eps}\times\R)\cap \Sphere^k_{2/\eps}}.
\]
\begin{lemma}\label{lem:Sh_eps}
For sufficiently small $\eps>0$ and $0\le i\le k$, we have
\[
b_i(S_{\eps})=\frac{1}{2}b_i(S^h_{\eps}).
\]
\end{lemma}
\begin{proof}
Note that $S_{\eps}$ is bounded by Proposition~\ref{prop:SvstildeS} and 
$S^h_{\eps}$ is the projection from the origin of the 
set~$S_{\eps}\times\{1\}\subset\R^k\times\{1\}$ onto the unit sphere 
in $\R^{k+1}$. Since $S_{\eps}$ is bounded, 
the projection does not intersect the equator and consists of two disjoint 
copies (each homeomorphic to the set $S_{\eps}$) 
in the upper and lower hemispheres.
\end{proof}
We now fix a sufficiently small $\eps>0$ and a 
family of polynomials that will be useful in what follows. 
Let 
\[
P=(\frac{2}{\eps})^2-\sum_{i=1}^{k+1} X_i^2.
\] 
By 
Proposition~\ref{prop:general} we can choose a family 
${\cH=\{H_1,\ldots,H_s\}}\subset\R[X_1,\ldots,X_{k+2}]$ of positive definite 
quadratic forms such that 
$\Zer(\cH_{J}\cup\{P^h\},\PP^{k+1}_{\C})$ 
is a non-singular complete intersection 
for every $J\subset\{1,\ldots,s\}$.

Let $\delta>0$ and let 
$\tilde{P}_i= (1-\delta)P^h_i + \delta \tilde{H}_i$, $1\le i\le s$, 
where $\tilde{H}_i=H_i(X_1,\ldots,X_{k+1},1)$. 
Note that $\tilde{H}_i$, $1\le i\le s$, is positive definite since $H_i$ is 
positive definite. 
Let $T_{\eps,\delta}\subset\R^{k+1}$ 
(resp., $\bar{T}_{\eps,\delta}\subset\R^{k+1}$) \
be the basic semi-algebraic set defined by 
${\tilde{P}_1>0,\ldots,\tilde{P}_s>0}$ (resp., ${\tilde{P}_1\ge 0,\ldots,\tilde{P}_s\ge 0}$) 
contained in $C_{\eps}$. 

Also, let 
\[
{\tilde{\cP}=\{\tilde{P}_1,\ldots,\tilde{P}_s\}}.
\]
\begin{lemma}\label{lem:T}
For all sufficiently small $0<\delta<\eps$ we have,
\begin{enumerate}
\item the homology groups of $S^h_{\eps}$ and $\bar{T}_{\eps,\delta}$ 
are isomorphic,
\item the homology groups of $T_{\eps,\delta}$ and $\bar{T}_{\eps,\delta}$ 
are isomorphic.
\end{enumerate}
\end{lemma}
\begin{proof}
For the first part note that $S_{\eps}$ and $\bar{T}_{\eps,\delta}$ 
have the same homotopy type 
using ~Lemma~16.17 in \cite{BPR06}.

The second part is clear since 
by choosing any slightly smaller $0 < \delta' < \delta$, 
we have a retraction from  $T_{\eps,\delta}$ to  $\bar{T}_{\eps,\delta'}$.
\end{proof}
Now, let $T^h_{\eps,\delta}\subset\R^{k+2}$ be the semi-algebraic set defined by 
${\tilde{P}^h_1>0,\ldots,\tilde{P}^h_s>0}$ contained in 
$\tilde{C}_{\eps}$, where 
\[
\tilde{C}_{\eps}=(\Ball^k_{1/\eps}\times\R^2)\cap (\Sphere^k_{2/\eps}\times\R)
\cap \Sphere^{k+1}_{4/\eps}.
\] 
Also, let 
\[
\tilde{\cP}^h=\{\tilde{P}^h_1,\ldots,\tilde{P}^h_s\}.
\]
\begin{lemma}\label{lem:sphvproj}
For all sufficiently small $0<\delta<\eps$ and $0\le i\le k$,
\begin{enumerate} 
\item\label{part1} $b_i(T_{\eps,\delta})=\frac{1}{2} b_i(T^h_{\eps,\delta})$,
\item\label{part2} for all $J\subset\{1,\ldots,s\}$, 
\[b_i\left(\Zer(\tilde{\cP}^h_J,\tilde{C}_{\eps})\right)=
2\cdot b_i\left(\Zer(\tilde{\cP}^h_J\cup\{P^h\},\PP^{k+1}_{\R})\right),
\]
\item\label{part3} for all $J\subset\{1,\ldots,s\}$,
$\Zer(\tilde{\cP}^h_J\cup\{P^h\},\PP^{k+1}_{\C})$ is a non-singular 
complete intersection.
\end{enumerate}
\end{lemma}
\begin{proof}
First, observe that 
$T^h_{\eps,\delta}\subset\tilde{C}_{\eps}$ 
is the projection from the origin of 
$T_{\eps,\delta}\times\{1\}\subset C_{\eps}\times\{1\}$ 
onto the sphere $\Sphere^{k+1}_{4/\eps}$ in $\R^{k+2}$. 
Note that 
$T^h_{\eps,\delta}$  
does not intersect the 
equator of the sphere $\Sphere^{k+1}_{4/\eps}$ in $\R^{k+2}$ (i.e. the set 
$\tilde{C}_{\eps} \cap \Zer(X_{k+2},\R^{k+2})$),
and consists of two disjoint copies (each 
homeomorphic to the set 
$T_{\eps,\delta}$) 
in the upper and lower hemisphere. 

For the second part, note that the set~$\Zer(\tilde{\cP}^h_J,\tilde{C}_{\eps})$ 
does not intersect the equator of the sphere $\Sphere^{k+1}_{4/\eps}$ in $\R^{k+2}$. 
Moreover, the two-fold covering
$\pi: \Sphere^{k+1}_{4/\eps} \rightarrow  \PP^{k+1}_{\R}$
(obtained by identifying antipodal points) 
restricts to a homeomorphism on the upper and lower hemisphere.

The third part follows from Proposition~\ref{prop:compl}.
\end{proof}
\begin{proposition}\label{prop:Tge0}
For all sufficiently small $0<\delta<\eps$ and for $0\le i\le k-1$, 
we have 
\[
b_i(T_{\eps,\delta}) \le  
\sum_{j=0}^{min\{s,k-i\}}{{s}\choose j}{{k+1}\choose {j}}2^{j}.
\]
\end{proposition}
\begin{proof}
By Lemma~\ref{lem:sphvproj}~(\ref{part1}) 
it suffices to prove the statement for the 
set~$T^h_{\eps,\delta}$. Note that 
$\Zer(\tilde{\cP}^h_J\cup\{P^h\},\PP^{k+1}_{\C})$ 
is a complete intersection for all $J\subset\{1,\ldots, s\}$ 
by Lemma~\ref{lem:sphvproj}~(\ref{part3}).

For $0\le i\le k-1$,
\begin{eqnarray*}
b_i(T^h_{\eps,\delta}) & \le & 
b_i\left(\Sphere^k_{4/\eps} \setminus \bigcup_{i=1}^s 
\Zer(\tilde{P}^h_i,\tilde{C}_{\eps})\right)\\
& \le & 1+ b_{k-1-i}\left(\bigcup_{i=1}^{s}
\Zer(\tilde{P}^h_i,\tilde{C}_{\eps})\right),
\end{eqnarray*}
where the 
first inequality is a consequence of the fact that,
$T^h_{\eps,\delta}$ is an open subset of 
${\Sphere^k_{4/\eps} \setminus \bigcup_{i=1}^s \Zer(\tilde{P}^h_i,\tilde{C}_{\eps})}$ 
and disconnected from its
complement in 
$\Sphere^k_{4/\eps} \setminus\bigcup_{i=1}^s \Zer(\tilde{P}^h_i,\tilde{C}_{\eps})$, 
and the 
last inequality follows from  Theorem~\ref{the:alex}~(Alexander Duality). 
It follows from 
Proposition~\ref{prop:spec}, Lemma~\ref{lem:sphvproj}~(\ref{part2}) 
and Theorem~\ref{the:smith}~(Smith inequality), 
that
\begin{eqnarray*}
b_i(T^h_{\eps,\delta}) & \le & 1+   \sum_{j=1}^{k-i}\sum_{|J|=j} 
b_{k-i-j}\left(\Zer(\tilde{\cP}^h_J,\tilde{C}_{\eps})\right)\\
& = & 1+ 2\cdot \sum_{j=1}^{k-i}\sum_{|J|=j} 
b_{k-i-j}\left(\Zer(\tilde{\cP}^h_J\cup\{P^h\},\PP^{k+1}_{\R})\right)\\
& \le & 1+ 2\cdot \sum_{j=1}^{k-i}\sum_{|J|=j} 
b_{k-i-j}\left(\Zer(\tilde{\cP}^h_J\cup\{P^h\},\PP^{k+1}_{\C})\right).
\end{eqnarray*}
Note that for $j\le s$ the number of possible $j$-ary intersections is equal to 
${{s}\choose j}$ and using Theorem~\ref{the:betti}, we conclude 
\begin{eqnarray*}
b_i(T^h_{\eps,\delta}) & \le & 1+ 2\cdot \sum_{j=1}^{min\{s,k-i\}}{{s}\choose j} b(j+1,k+1)\\ 
& \le &  
2\cdot \sum_{j=0}^{min\{s,k-i\}}{{s}\choose j}{{k+1}\choose {j}}2^{j}.
\end{eqnarray*}
The claim follows since 
$b_i(T_{\eps,\delta})=\frac{1}{2}\cdot b_i(T^h_{\eps,\delta})$.
\end{proof}
%
\pagebreak[2]
We are now in a position to prove our main result.
\begin{proof}[Proof of Theorem~\ref{the:Pgre0}]
For all sufficiently small $0<\delta<\eps$ we have by Lemma~\ref{lem:T} 
that the homology groups of 
$S^h_{\eps}$ and $T_{\eps,\delta}$ are isomorphic. Moreover, 
for $0\le i\le k-1$, $b_i(S)=\frac{1}{2}b_i(S^h_{\eps})$ by 
Proposition~\ref{prop:SvstildeS} and Lemma~\ref{lem:Sh_eps}. 
Hence, the first part 
follows from Proposition~\ref{prop:Tge0}. 

The second part follows from an easy computation.
\end{proof}
Finally, we prove Corollary~\ref{cor:bsum}.
\begin{proof}[Proof of Corollary~\ref{cor:bsum}]
Follows by applying the bound of Theorem~\ref{the:Pgre0} 
to each $b_i(S)$, ${0\le i\le k-1}$.
\end{proof}

\section{Conclusion and Open Problems}
\label{sec:open}
In this paper we have improved the upper bound proved by Barvinok on
the Betti numbers of semi-algebraic sets in $\R^k$ defined by 
$s \leq \frac{k}{2}$ 
quadratic inequalities. 
The new bound is of the form $(O(k/s))^s$ 
improving the previous bound
of $k^{O(s)}$ due to Barvinok. Using the fact that
a complex non-singular complete intersection in $\C^k$ defined by $s$
quadratic equations can be viewed as a real semi-algebraic set in $\R^{2k}$ 
defined $2s$ quadratic equations, it follows that the best
bound on the sum of the Betti numbers of semi-algebraic sets defined by $s$ quadratic
inequalities in $\R^k$ cannot be better than $k^{O(s)}$. We conjecture that
the exponent, $s$, in our bound is in fact optimal and an interesting
open problem is to construct an example which meets our bound. 

Another interesting problem in this context is to obtain a tighter bound
on the number of connected components (that is on $b_0(S)$) for $S \subset
\R^k$ defined by $s \leq k$ quadratic inequalities. It can be easily seen
from the example of $S \subset \R^k$ defined by,
\[ X_1(X_1 - 1) \geq 0, \ldots, X_s(X_s - 1) \geq 0, \]
that $b_0(S)$ can be as large as $2^s$. However, we know of no examples
where $b_0(S)$ is as large as $k^{\Omega(s)}$. 
Note that the Betti numbers of a non-singular complex complete intersection is 
concentrated in the ``middle'' dimension. Consequently,
the Smith inequality  gives bounds only on the sum of the
Betti numbers of the corresponding real varieties.
Because of this the method of the proof used in this paper
has the drawback that it gives no way of proving better 
bounds on the individual (say the lowest or the highest) Betti numbers. 
%

%
\end{document}